\documentclass[11pt]{article}
\usepackage[english]{babel}
\usepackage{amsmath,amssymb, amsthm}
\usepackage{graphicx}
\usepackage{tikz}
\usepackage{url}
\usepackage{lineno}
\usepackage[a4paper]{geometry}
\newtheorem{theorem}{Theorem}

\newtheorem{lemma}[theorem]{Lemma}
\newtheorem{remark}[theorem]{Remark}

\def\RR{{\mathbb{R}}}
\def\ZZ{{\mathbb{Z}}}

\def\phi{\varphi}
  \newcommand{\N}{\ensuremath{\mathbb{N}}}

\newcommand{\R}{\ensuremath{\mathbb{R}}}

\begin{document}
\begin{center}
\Large  {\small\date{today} } \smallskip
 \bf New methods for quasi-interpolation approximations:
  resolution of odd-degree singularities
\end{center}
\bigskip\bigskip
\begin{center} \it Martin Buhmann, Justus-Liebig University,
  Mathematics, 35392 Giessen, and \\ \smallskip
  Janin J\"ager, 
Katholische Universit\"at Eichst\"att-Ingolstadt, 85049 Ingolstadt,
Germany,\\ \smallskip 
Joaqu\'in J\'odar,
Department of Mathematics, University of Ja\'en, and\\ \smallskip
Miguel L. Rodr\'iguez,
Department of Mathematics, University of Granada, Spain.

\end{center}
\bigskip\bigskip
{\small {\it Abstract: \/} In this paper, we study functional approximations
where we choose the so-called radial basis function method and more
specifically, 
quasi-interpolation. From the various available approaches to the latter, we
form new quasi-Lagrange functions when the orders of the singularities of the radial
function's Fourier transforms at zero do not match the parity of the dimension of the space,
and therefore new expansions and coefficients are needed to overcome
this problem. We develop explicit constructions of infinite Fourier
expansions that provide these coefficients and make an extensive
comparison of the approximation qualities and - with a particular
focus - polynomial precision and uniform approximation order of
the various formulae. One of the interesting observations concerns the
link between algebraic conditions of expansion coefficients and
analytic properties of localness and convergence.}

\date{\today}  
 \bigskip\bigskip\bigskip
 \title{Introduction}
 \section{Introduction}\label{rem}
In this paper, we study functional approximations in one or more
variables to approximands that are at a minimum continuous. A variety
of approximation methods are available by splines, multivariable
polynomials, trigonometric polynomials etc., but here we choose the
so-called radial basis function method.

Among the approximants that can be formed from the linear spaces
spanned by shifts $\phi(\|\cdot-x_j\|)$, where the norm is usually
Euclidean and $\phi:\RR_+\to\RR$ is the so-called radial basis
function, are mostly interpolants (see \cite{rbf}) and
quasi-interpolants  (see \cite{qi}). Those in
turn can be formed for finitely many scattered data or gridded data,
often infinitely many.

We will study quasi-interpolants in this article because they are
relatively simple to form and have excellent convergence
properties. Those originate from the spaces spanned by shifts
$\psi(\cdot-x_j)$ containing polynomials of some low  total degree
{\it and\/} the $\psi$ themselves coming from the span of shifts
$\phi(\|\cdot-x_j\|)$ decaying quickly. Also, collocation if often
explicitly not desired. Those two properties together
avail us to establish convergence theorems if the approximands $f$ are
smooth enough, the basis ingredients coming from local Taylor
expansions of course.

All this rests on the ``quasi-Lagrange functions'' satisfying the
famous Strang and Fix conditions, which we cite in an appropriate form for the convenience of the reader at the end of the introduction. A basic tool is that of Fourier transforms, in one or more dimensions, of a function $f$,
$$
\hat{f}(\xi)=\int_{\mathbb{R}^n} f(x) e^{-i \xi \cdot x} dx,\quad \xi\in\mathbb{R}^n.
$$

Since the $\psi$ are to be linear combinations of the $\phi(\|\cdot-x_j\|)$ satisfying these conditions depends on orders of
singularities of the generalised Fourier transforms of the radial
basis functions $\phi$. In order to resolve those singularities at the
origin, the coefficients forming the quasi-Lagrange functions from the
$\phi$s are designed such that there is a high-order contact between
the functions at zero.

This works very well if the parities of the singularities and the
orders of the zero of the trigonometric polynomials in Fourier space
that are formed from the aforementioned coefficients match and are
even, but if the orders of the singularities are odd, extra work has
to be done. Roughly speaking, it is no longer possible to have
trigonometric {\it polynomials\/} coming from the said coefficients of
$\psi$ as a linear combination of shifts $\phi(\|\cdot-x_j\|)$, but
they end up in infinite series in order to have the required high
order contact absorbing the singularities of different parities.

We aim to undertake an analysis of this phenomenon, our prime examples
being multiquadrics (which ought to be used in odd dimensions so that
the order of the singularity at the origin is even) and thin-plate
splines (which ought to be used in even dimensions so that
the order of the singularity at the origin is even). And we intend to
employ these radial functions just in the spaces of ``wrong''
dimensions, i.e., even and odd respectively and develop and compare
different methods to resolve the mentioned problems (cardinal
functions, infinite expansions for quasi-Lagrange functions and
others). The specific examples we study are the thin-plate spline in dimension one and then compare to the generalised multiquadric in one dimension, which is in the 'right' dimension. For the latter we give an explicit expression of the quasi-Lagrange function. We put particular emphasis not only on comparisons but on
the order of decay (localness) of the $\psi$s and of course on the  
polynomial precision of the generated vector spaces.  We shall need the following famous result of Strang and Fix several times. (See, among many other possible references, \cite[Theorem
  2.2]{qi}). 
\begin{theorem}\label{SF}[Strang and Fix conditions]

Let $\psi:\mathbb{R}^n \rightarrow
\mathbb{R}$ be a continuous function such that
\begin{enumerate}
\item there exists a positive $\ell$ such that for some nonnegative integer $m$, when
  $\|x\|\rightarrow \infty$, $|\psi(x)|={O}(\Vert
  x\Vert^{-n-m-\ell})$, which immediately implies $m$-fold
  differentiability of the Fourier transform,
\item $D^{\alpha} \hat{\psi}(0)=0$, $\forall \alpha \in \mathbb{Z}_+^n$, $1\leq \vert \alpha\vert\leq m$, and $\hat{\psi}(0)=1$,
\item $ D^{\alpha}\hat{\psi}(2\pi j)=0,\ \forall j \in \mathbb{Z}^n\setminus\{0\}$ and  $\forall \alpha \in \mathbb{Z}_+^n$ with $ \vert\alpha\vert \leq m$. 
\end{enumerate}

Then the quasi-interpolant
\begin{equation}
Q_hf(x)=\sum_{j\in\mathbb{Z}^n}f(hj)\psi(x/h-j),\qquad x\in\RR^n,
\end{equation}
is well-defined and exact on $\mathbb{P}_m$. The approximation error can be estimated by
$$\Vert Q_hf-f\Vert_{\infty}=\begin{cases}{O}(h^{m+\ell}), & \text{when } 0<\ell<1,\\
{O}(h^{m+1}\log (1/h)), & \text{when } \ell=1,\\
{O}(h^{m+1}), & \text{when } \ell>1,\\
\end{cases}$$
for $h\rightarrow 0$ and a bounded function $f\in C^{m+1}(\RR^n)$ with bounded derivatives.
\end{theorem}
Here and below  $\mathbb{P}_m$ denotes the space of polynomials of
degree at most $m$.

\section{The thin-plate spline in one-dimension using infinitely many coefficients}\label{sect2}

Along this section and as outlined in the introduction, we consider the radial basis function:
$\varphi(r)=r^2\log (r)$. Therefore the relevant generalised Fourier
transform is (see \cite[Chapter 4.6]{jon}):
$$\hat{\varphi}(r)=-2\pi\left(\left(1+\frac12-\gamma\right)\delta^{''}(r)- r^{-3}\right),$$ 
where, as usual,  $\delta$ denotes the Delta-distribution  and  $\gamma$ is Euler's constant ($\gamma \equiv 0.577216$).


At this step, we define $\mu$ as the real number (if it exists) satisfying $\hat{\varphi}(r) \doteq r^{-\mu}+O(r^{-\mu+1})$, $r\rightarrow 0^+.$
Here, the $\doteq$ means equality up to a nonzero constant multiple.

We have $\mu=3$  in this particular case.  The quasi-Lagrange functions will take the form
$$\psi(x)=\sum_{j\in \mathbb{Z}} \lambda_j \phi(|x-j|),\qquad
x\in\mathbb{R}.$$
We study the resulting schemes for different choices of the
$\lambda$ coefficients in this ``wrong'' (odd)  dimension, where we
cannot achieve the Strang-Fix conditions using only finitely many
coefficients because the trigonometric expansions with odd-order zeros
at the origin (resolving odd order singularities) will always be
infinite unlike the even power singularities (we can use powers and
tensor-products of $1-\cos\xi$). 
\subsection{Cardinal interpolation, infinite expansions from shifts of
the radial function}

An Ansatz that always works in forming of Lagrange functions (no
longer quasi-Lagrange functions) from equally spaced shifts
$$\psi(x)=\chi(x)=\sum_{j\in \mathbb{Z}} \lambda_j \phi(|x-j|),\qquad
x\in\mathbb{R},$$
which satisfy $\psi(k)=\chi(k)=\delta_{0k}$ for all integers $k$,
where we use the standard notation $\chi$ for the cardinal functions. 

The coefficients 
$\lambda_j$ come from the Fourier expansion within the Wiener algebra
of the reciprocal of the so-called symbol,
$\sigma(\vartheta)=\sum_{\ell \in \ZZ} \hat{\phi}(\vartheta+2\pi \ell)$, so that
$$\lambda_j= \displaystyle \frac1{2\pi}\int_{\mathbb{T}}\dfrac{e^{i \vartheta
    j}}{\sigma(\vartheta)} {\rm d}\vartheta,\qquad j\in\ZZ,$$ and
$$\psi(x)= \displaystyle \sum_{j\in \mathbb{Z}} \lambda_j \varphi(|x-j|),\qquad x\in\RR.$$
The localness or asymptotic decay of $\psi(x)$ is identified as
$|\psi(x)|=O((1+|x|)^{-4}))$ for all (in particular large in modulus) $x$. For a
suitable result see \cite[Theorem 4.3]{rbf}, and here $\mu=3$. We use
$\int_{\mathbb{T}}$
for $\int_{-\pi}^\pi$ throughout.

As we stated in the introduction, polynomial precision 
$$ \sum_{k\in\mathbb{Z}} p(k)\psi(\cdot-k)\equiv p$$
for some (usually, low-degree) polynomials $p$ is crucial. In this
case we get the polynomial precision/reproduction:
$\mathbb{P}_2$-reproduction because $\mu=3$. (See \cite[Theorem
  4.4]{rbf}).

A by now standard theorem then delivers a uniform approximation error for
suitably smooth approximands $f\in C^{4}(\mathbb{R})$ with bounded
derivatives
$$ \Bigl\| f-  \sum_{k\in\mathbb{Z}}
f(kh)\psi(h^{-1}\cdot-k))\Bigr\|_\infty=O(h^3 |\log h|)$$ for $h\to
0$. (See \cite[Theorem 4.6]{rbf} or the convergence results in \cite{BDai}.) We have to apply the remark in the
proof, and not Theorem 4.6  itself, because $\mu$ is an integer or we
can apply Theorem \ref{SF} with  $m=2, \ell=1, n=1$.

While cardinal interpolation always works (note the absence of demands
for certain parities of dimensions and orders of singularities)
because  no polynomials of trigonometric type but infinite expansions
are used, we now turn to ``genuine'' quasi-interpolations where no
cardinal conditions are demanded.

\subsection{Quasi-interpolation without cardinality conditions}\label{subsect22}
We now wish to go away from the well-known cardinal function approach
and use straight quasi-interpolation instead. That brings us to the
problem, when the parity of the radial function's generalised Fourier
transform's singularity at zero is odd, we can no longer form a trigonometric {\it
  polynomial\/} $q$, say, that matches the degree of the said
singularity. This comes from the fact that only even powers of
trigonometric expansions have finitely many coefficients when written
as Fourier series expansions.  So we need to use expansions of
periodic series with
infinitely many coefficients $\lambda_j$ coming from the Fourier
coefficients of $(2-2\cos x)^{3/2}$ or for instance $|\sin
x|^3$. These are of course by no means unique, but mere
examples. Therefore we end up in expressions
$$\psi(x)= \displaystyle \sum_{j\in \mathbb{Z}} \lambda_j
\varphi(|x-j|),\qquad x\in\RR.$$

We like to apply the famous Strang and Fix conditions to check the
polynomial precision's degree of quasi-interpolation with thin-plate
splines in one dimension
$$ \sum_{k\in\mathbb{Z}} p(k)\psi(\cdot-k)\equiv p.$$
They depend on the properties of the (classical) Fourier transforms of $\psi$.
By straightforward computations (and dividing by $2\pi$ in order to normalise) it can be shown that
$$\hat{\psi} (0)=1, \qquad  \dfrac{d\, \hat{\psi}}{d \xi}(0)=0, \qquad \dfrac{d^{k}\hat{\psi}}{d \xi^k}(2\pi j)=0, \: \forall j \in \mathbb{Z}\setminus \{0\}, \: k=0, 1, \quad \dfrac{d^2\, \hat{\psi}(\xi)}{d \xi^2}(0)\neq 0. $$ 
As the Strang-Fix condition for first degree derivative is satisfied
but the second derivative is not satisfied at $\xi=0$, we have
$\mathbb{P}_1-$reproduction at a maximum (but not any higher). We arrive at
\begin{theorem} Let $\phi$ be the thin-plate spline radial basis
  function. With $n=1$ and the Fourier coefficients of  $(2-2\cos
  x)^{3/2}$ being the $\lambda_j$, the quasi-interpolation 
  $$ \sum_{k\in\mathbb{Z}} p(k)\psi(\cdot-k)\equiv p$$
  is exact for linear polynomials $p$.
  \end{theorem}
The next question is the decay of the quasi-Lagrange functions $\psi$:
We expect $|\psi(x)|=O((1+|x|)^{-4})$  (see Theorem \ref{th1} for
$\underline n=3$). This is one order better than the routinely
required third order decay which would suffice for absolute
convergence of the quasi-interpolant when at most linearly growing
approximands (in particular linear polynomials) are inserted. 

This leads us to the question of approximation error; a routine result
gives us from the identified polynomial precision and the order of
decay of the quasi-Lagrange functions $O(h^2 )$ for $h\to 0$. (See
Theorem \ref{SF} with  $m=1, \ell=2, n=1$.) Notice the absence
of the logarithmic term due to the one order faster decay of the quasi-Lagrange function.

\subsection{An intermediate formulation of $\psi$ and its Fourier transform}\label{subsect23}
Another scheme is the separation of the Fourier transform's
singularity into two factors: one that resolves the odd singularities
degree separately and leaving an even negative power, and then using
the classical approach for the high even order contact at the origin.

So we begin in trying to improve the polynomial reproduction using the
above scheme and the situation at zero. We set therefore $$\hat{\psi}(\xi)=P(\xi) |\sin(\xi)| \hat{\varphi}(|\xi|),\qquad\xi\in\RR,$$
being $P(\xi)= \displaystyle \sum_{k=-N}^N \mu_k e^{ik\xi}$ a suitable
trigonometric polynomial. In summary, we set this new scheme in this way: 

We have to fix the coefficients of the quasi-Lagrange functions. This
always begins with setting the coefficients $\lambda_k$ to be the
(infinitely many) Fourier coefficients of the expansion of
$$P(\xi)|\sin (\xi)|,\qquad -\pi\leq \xi\leq\pi,$$ -- this
  is by no means unique,  we could use for instance $$P(\xi)(1-\cos
  (\xi))^{1/2},\qquad -\pi\leq \xi\leq\pi.$$
  Many other choices of roots of trigonometric functions are possible.
Depending on those coefficients, especially which trigonometric expansions they
form and which orders their zeros have, we will arrive at a reproduction of polynomials $\mathbb{P}_2-$precision (see below).

In order to compute the terms and derivatives that will serve to
verify up to which order the Strang-Fix  conditions hold at zero, we will use 
$$\hat{\psi}(\xi)=2\pi\dfrac{P(\xi) |\sin(\xi)|}{|\xi|^3}\equiv
2\pi\dfrac{P(\xi) \sin(\xi)}{\xi^3} $$ 
as $\xi$ goes to 0 because the singular term in $\hat{\varphi}(r)$ is a constant multiple of $r^{-3} $ at the
origin.  In order to satisfy the conditions at the zero  for
$\mathbb{P}_2$-reproduction   we demand according to the Strang and
Fix approach
\begin{equation}\label{strang}
\hat{\psi}(\xi)=1+ O(\xi^3).
\end{equation}
Therefore, close to  the origin, we can expand
$$2\pi\dfrac{P(\xi) \sin(\xi)}{\xi^3}=2\pi \displaystyle
\sum_{k=-N}^N \mu_k e^{ik\xi}\left (\xi - \dfrac{\xi^3}{3!}+
\dfrac{\xi^5}{5!}- \dfrac{\xi^7}{7!} \pm \cdots \right)\times \xi^{-3}$$
which is
$$2\pi\displaystyle \sum_{k=-N}^N  \sum_{j=0}^{\infty} \mu_k \dfrac{(ik\xi)^j}{j!}\left (\xi - \dfrac{\xi^3}{3!}+ \dfrac{\xi^5}{5!}- \dfrac{\xi^7}{7!} \pm \cdots \right)\times \xi^{-3}.$$
Now, by imposing the  condition (\ref{strang}) we have
 $$1+ O(\xi^3)=
 2\pi\displaystyle \sum_{k=-N}^N  \sum_{j=0}^{\infty} \mu_k \dfrac{(ik\xi)^j}{j!}\left (\xi^{-2} - \dfrac{1}{3!}+ \dfrac{\xi^2}{5!}- \dfrac{\xi^4}{7!} \pm \cdots \right).$$
The above equation will give  conditions on the  $\mu_k$s  (they will not be unique, of course). Specifically, for $N=2$ we obtain \begin{equation*}
\displaystyle \sum_{k=-2}^2 \mu_k=0,\:
\displaystyle \sum_{k=-2}^2 k \mu_k=0, \: 
- \pi\displaystyle \sum_{k=-2}^2 \left(k^2+ \dfrac13\right) \mu_k
= 1, \:
\end{equation*}
and further
\begin{equation*}
- \dfrac16\displaystyle \sum_{k=-2}^2 \left(k+ k^3\right) \mu_k = 0,\:
\dfrac{1}{12}\displaystyle \sum_{k=-2}^2 \left(\dfrac{1}{10}+ k^2+ \dfrac{k^4}{2}\right) \mu_k = 0.
\end{equation*}
This can be formulated equivalently as 
\begin{equation} \label{mu}
\displaystyle \sum_{k=-2}^2 \mu_k=
\displaystyle \sum_{k=-2}^2  k \mu_k=0, \quad
\displaystyle \sum_{k=-2}^2 k^2 \mu_k = -\dfrac{1}{\pi}, \quad
\displaystyle \sum_{k=-2}^2 k^3 \mu_k = 0, \quad
\displaystyle \sum_{k=-2}^2 k^4 \mu_k = \dfrac{2}{\pi}.
\end{equation}
We will then choose
 $$\mu_{-2}= \dfrac{1}{8\pi}, \:  \mu_{-1} = -\dfrac{1}{\pi}, \:  \mu_0 =
\dfrac{7}{4\pi}, \: \mu_1 =  -\dfrac{1}{\pi}, \: \mu_2= \dfrac{1}{8\pi} $$
or
$$P(\xi)=  \dfrac{1}{8\pi} e^{-2 i \xi} -\dfrac{1}{\pi} e^{-i \xi}+ \dfrac{7}{4\pi} -\dfrac{1}{\pi} e^{i \xi}+\dfrac{1}{8\pi} e^{2 i \xi},$$ i.e.
 $$P(\xi)= \dfrac{7}{4\pi}-\dfrac{2}{\pi} \cos(\xi) + \dfrac{1}{4\pi}\cos(2\xi).$$
We can easily verify  that the upper bound on the right-hand side of (\ref{strang}) holds for $\hat{\psi}(\xi)=P(\xi)|\sin (\xi)|\hat{\varphi}(|\xi|)$.

Now,  let us study the behaviour of the derivatives $\hat{\psi}^{(\ell)}
(2\pi j)$, $j\in\mathbb{Z}\setminus \{0\}$. The purpose of this is of
course checking the Strang and Fix conditions:
\begin{enumerate}
\item For $\ell=0$ we have
\begin{equation*}
\hat{\psi}(2\pi j)=0,\quad j\in\mathbb{Z}\setminus \{0\},
\end{equation*}
because of the $|\sin (\xi)|$-term and the continuity of $\hat{\varphi}(\xi)$
away from the origin. 
\item For $\ell=1$
\begin{equation}\label{s1}
\dfrac{d\hat{\psi}}{d \xi}(\xi)=P'(\xi)|\sin (\xi)|\hat{\varphi}(\xi) +P(\xi) \dfrac{\sin(\xi) \cos(\xi)}{|\sin (\xi)|}\hat{\varphi}(\xi)+ P(\xi)|\sin (\xi)|\hat{\varphi}'(\xi)
\end{equation} 
which vanishes for all $2\pi j$, $j \in \mathbb{Z}\setminus \{0\}$: the first and the third term clearly vanish;  the second term of the right-hand part of  (\ref{s1}) vanishes due to $\displaystyle \sum_{k=-2}^2 \mu_k=0$.

\item And for $\ell=2$ the derivatives of the first and third terms of
  (\ref{s1})  vanish at $2\pi j$ with $j \in \mathbb{Z}\setminus
  \{0\}$ (we  now have to use  that $\displaystyle \sum_{k=-2}^2k
  \mu_k=0$).  In fact, problems could come from the derivative of the
  second term of (\ref{s1}) which we will therefore have to compute
  explicitly. It is
\begin{equation*}
\begin{split}
&P'(\xi) \dfrac{\sin(\xi) \cos(\xi)}{|\sin (\xi)|}\hat{\varphi}(\xi) +P(\xi) \dfrac{\cos^2(\xi)}{|\sin (\xi)|}\hat{\varphi}(\xi)
-P(\xi) |\sin (\xi)|\hat{\varphi}(\xi)\\
&
+P(\xi) \dfrac{\sin(\xi) \cos(\xi)}{|\sin (\xi)|}\hat{\varphi}'(\xi)
-P(\xi) \dfrac{\cos^2(\xi)}{|\sin (\xi)|}\hat{\varphi}(\xi)\\
&=P'(\xi) \dfrac{\sin(\xi) \cos(\xi)}{|\sin (\xi)|}\hat{\varphi}(\xi)
-P(\xi) |\sin (\xi)|\hat{\varphi}(\xi)
+P(\xi) \dfrac{\sin(\xi) \cos(\xi)}{|\sin (\xi)|}\hat{\varphi}'(\xi).
\end{split}
\end{equation*}
In summary, the second order derivatives of $\hat{\psi}$ vanish at the points $2\pi j$, $j\in\mathbb{Z}\setminus\{0\}$ and we have all the requirements for getting  $\mathbb{P}_2-$reproduction.
\end{enumerate}
We sum our findings up in the following result.
\begin{theorem} Let $\phi$ be the thin-plate spline radial basis
  function. With $n=1$ and the Fourier coefficients of
  $P(\cdot)\times|\sin(\cdot)|$, $P$ as above,
  being the quasi-Lagrange function's coefficients $\lambda_j$, the quasi-interpolation 
  $$ \sum_{k\in\mathbb{Z}} p(k)\psi(\cdot-k)\equiv p$$
  is exact for quadratic polynomials $p$.
  \end{theorem}
\begin{remark} At the points $x=2\pi j, j \in \mathbb{Z}\setminus
  \{0\}$ the third derivative of $\hat{\psi}$ has a jump discontinuity,
  so we will not be able to satisfy high enough degree Strang-Fix conditions that we
  could successfully impose in order to arrive at $\mathbb{P}_3-$reproduction.
\end{remark}
We note that a convenient way to show that the second order
Strang-and-Fix conditions, and only those degrees, are satisfied at the $\xi=2\pi j, j \in \mathbb{Z}\setminus
  \{0\}$ (while they would hold at zero up to order three, although
  that does not help) is to notice the following:
   if we define
  $g(\xi)=\frac12\sin|\xi|+\frac12\sin|\xi-\pi|$, then
  $$|\sin \xi| = \sum_{k=-\infty}^\infty g(\xi-k\pi),\qquad
  \xi\in\mathbb{R}.$$
  Therefore, using $P$ as above and expanding about the origin we get
  near zero
  \begin{align*}
P(\xi)|\sin \xi|\hat\varphi(|\xi|)&= 2 \pi P(\xi)|\sin \xi| \times |\xi|^{-3}
=(\sin \xi) 2 \pi \xi^{-3} \Bigl( \frac{\xi^2}{2 \pi } +\frac{\xi^4}{12 \pi }  -\frac{7\xi^6}{360 \pi } +\cdots\Bigr)
\\&
= \Bigl( \xi- \frac{\xi^3}{6} + \cdots  \Bigr) \times \Bigl( \frac{1}{\xi} +\frac{\xi}{6}-\frac{7\xi^3}{180}   +\cdots\Bigr)=
  \end{align*}
  which is
  $$1-\frac{7}{120}\xi^4 + O(\xi^6)$$ near zero (so third order SF-conditions possible)
  but since
  $$P(\xi)|\sin \xi|\hat\varphi(|\xi|)=|\xi|^{-3} \left( |\xi-2\pi j|^{3} + O( |\xi-2\pi j|^{7})\right)$$ near $\xi=2\pi j, j \in \mathbb{Z}\setminus
  \{0\}$, only second order Strang-and-Fix conditions are satisfied. The
  mentioned third degree derivative discontinuity can be seen as well.

  In fact, using the above form of $|\sin|$ we can compute its
  generalised Fourier transform which could then be used to compute
  some $\psi$ explicitly.

The details of the computation are given in Appendix A.1 Lemma \ref{sin}, where we show that 
  $${\cal F}^{-1}|\sin|(x)=\frac1{\pi}\times\frac{1+\exp(ix\pi)}{1-x^2}\times{\cal D}_2(x).$$
  Here ${\cal D}_2$ is the Dirac comb
  $${\cal D}_2=\sum_{k=-\infty}^\infty \delta(\cdot-2k).$$
  We always have to study the decay of $\psi$, because the resolution of
  the singularities at the origin will not deliver the desired
  polynomial precision of $$ \sum_{k\in\mathbb{Z}}
  p(k)\psi(\cdot-k)\equiv p$$
  unless the series above converge absolutely. For this we will need
  at least an asymptotic decay of $O((1+|x|)^{-3-\varepsilon})$ for
  the quasi-Lagrange function in order to get the summability of the
  series for the aforementioned quadratic polynomial reproduction. The asymptotic decay is
  fairly easily established (see, for instance, \cite{rbf} or
\cite{qi}) by exploiting the differentiability
  properties of the quasi-Lagrange function's Fourier transform. At
  even order multiples of $\pi$ this $\hat\psi$ is infinitely smooth,
  but at odd multiples of $\pi$, we observe the following
  behaviours. We get about $\xi=(2j+1)\pi$ (but not about $\xi=2j\pi$)
$$\lim_{\xi \to -\pi^-} \hat{\psi}'(\xi)=-\dfrac{8}{\pi^3}, \quad \lim_{\xi \to -\pi^+} \hat{\psi}'(\xi)=\dfrac{8}{\pi^3}.$$
As $\hat{\psi}(\xi) \notin C^1(\mathbb{R})$, then the maximum decay we can obtain for $\psi(x)$ will be $O((1+|x|)^{-2})$ according to our previous remark.

Finally we note that the approximation order (with the established decay and
  \cite[Theorem 2.2]{qi} with $m=0$) will be $O(h\, |\log h|)$.

%
%
%
%

\begin{remark}
So far, we have only considered pointwise function evaluation as a means to put
the approximand's information into the quasi-interpolant. It is
entirely possible to formulate the latter with different linear
functionals applied to the approximation:
 $$Q f(x)= \displaystyle \sum_{j\in \mathbb{Z}} \lambda_j(f) \psi(x-j)$$ where $\lambda_j(f)$ are suitable functionals. For example, they could be related to the  polynomial $P(x)$, i.e. 
$$\lambda_j(f)=\sum_{k=-M}^{M} \mu_k f(j-k) \quad {\rm or} \quad
\displaystyle \lambda_j(f)=\sum_{k=0}^{M}\mu_j f^{(k)}(j)$$ being
$\mu_k$ the coefficients of $P(x)$ in order to get polynomial
reproduction; we may also take local integrals.
  
That Ansatz would amount to a preconditioning of the
approximand before it is fed into the quasi-interpolation procedure.

As soon as the decay of $\psi$ in absolute terms is not
sufficient for the absolute convergence (summability) of the series in
$$ \sum_{k\in\mathbb{Z}} p(k)\psi(\cdot-k)\equiv p,$$
the functionals must help in order to improve the decay of the
function $\psi$ because   we need the summability of the series so
that the polynomial precision is formed in a well-defined
way. However, we do not provide any further detailed analysis of this
aspect in this work (but see \cite{qi}).
\end{remark}

\section{Finitely many coefficients for the quasi-Lagrange functions}\label{sect3}
\subsection{The classical scheme with the ``natural degree'' radial
  basis function in $\mathbb{R}$}\label{2s3}
In odd dimension the ``natural degree'' (i.e., even order singularity
at the origin) radial basis function using
multiquadrics and their ilk will be  $\varphi(r)=(r^2+ c^2)^{3/2}$ instead $\varphi(r)=r^2\log (r)$. 
In this case we  have the following scheme.

The radial basis function is in the first place $\varphi(r)=(r^2+
c^2)^{3/2}$ which we shall call the generalised multiquadric function.  Its Fourier
transform has the property that near zero we have that
$\hat{\varphi}(r)\doteq r^{-4}$, as required (even order), and therefore in our notation above
$\mu=4$. 

The quasi-Lagrange function's now available {\it finite number\/} of
coefficients  $\lambda_j$  are the Fourier coefficients of $(1-\cos
(\xi))^{2}$. The latter is a trigonometric polynomial with an even
order zero at the origin. Thus,  the
$\lambda_j$ will have finite support with respect to their indices. 
One way of choosing the Fourier transform is, as a consequence, $\hat{\psi}(\xi)=(1-\cos
(\xi))^{2} \hat{\varphi}(|\xi|)$.

It is interesting to verify the polynomial precision results for the
cases in the classical way of analysis.
Taking into account that the (distributional) Fourier transform, $\hat{\Phi}(\xi)$, of the generalised multiquadric $\Phi(x)=(c^2+\|x\|^2)^{\beta}$, $x\in\mathbb{R}^n$, $c>0$, $\beta\in\mathbb{R}\setminus\mathbb{N}_0$ is
$$
\hat{\Phi}(\xi)=(2\pi)^{n/2}\frac{2^{1+\beta}}{\Gamma(-\beta)} \left(
\frac{c}{\|\xi\|} \right)^{\beta+\frac{n}{2}}\times K_{\beta+\frac{n}{2}} (c\|\xi\|),\quad \xi \neq 0,
$$
according to \cite{jon} we have, in our case ($n=1, \beta=3/2$),
$$\hat{\varphi}(\xi)=(2\pi)^{1/2} \frac{2^{5/2}}{\Gamma (-3/2)}
\frac{c^2}{\xi^{2}} K_2(c|\xi|).$$
Furthermore, from  \cite{olv} we read
\begin{align*}\label{Olver}
\frac{c^s}{\|\xi\|^{s}} K_s(c\|\xi\|) =& 2^{s-1}\frac{1}{\|\xi\|^{2s}}
\sum_{k=0}^{s-1} \dfrac{ (s-k-1)!}{k!(-4)^k}(c\|\xi\|)^{2k}+{}\\
&{}
-\left(\frac{-c^2}{2}\right)^s   \log \|\xi\|    \sum_{k=0}^{\infty}
\frac{(c\|\xi\|)^{2k}}{4^kk!\Gamma(s+ k+1)}  +{}
\\&{}-
\left(\frac{-c^2}{2}\right)^s\log c   \sum\limits_{k=0}^{\infty} \frac{(c\|\xi\|)^{2k}}{4^kk!\Gamma(s+ k+1)} +{}\\&+{} \left(\frac{-c^2}{2}\right)^s \sum_{k=0}^{\infty} \left(\frac{\log 2}{4^kk!\Gamma(s+
  k+1)}+\frac12\dfrac{\Psi(k+1)+\Psi(s+k+1)}{4^k(s+k)!k!}\right)\times{}\\
&{}\times(c\|\xi\|)^{2k},
\end{align*}
where $\displaystyle\Psi(z)=\frac{\Gamma '(z)}{\Gamma (z)}$ is the Digamma function, and $\Gamma$ is the Gamma function.
 We summarise it in the following result.
\begin{theorem}
The quasi-Lagrange function satisfies the bound
$|\psi(x)|=O((1+|x|)^{-5})$.
\end{theorem}
This localness, i.e.,
decay of $\psi$, is  identified by differentiating its Fourier
transform and because 
 $\hat{\psi}(\xi)$ can be expanded, about $\xi=0$, as
$$\tilde{a}+a\xi^2+b\xi^4\log|\xi|+\cdots, \quad \tilde{a},a,b\in\mathbb{R},$$
where the dots mean higher order terms (higher powers of $\xi$ and/or logarithms) in the
expansion. Indeed, this
comes from the previous expansion of $\frac{c^2}{\|\xi\|^{2}}
K_2(c\|\xi\|)$ and the expansion of $(1-\cos (\xi))^{2}$ given by
$\frac{\xi^4}{4}$, around $\xi=0$. Therefore, by applying the (inverse) Fourier
transform we would obtain the above-mentioned decay (see \cite{jon},
for instance). 

In order to prove that the  polynomial reproduction is $\mathbb{P}_1$,
we make some explicit computations. 
First we derive the suitable normalisation constant, by noting that
\begin{align*}\hat{\varphi}(|\xi|)&=(2\pi)^{1/2} \frac{2^{5/2}}{\Gamma (-3/2)}
\frac{c^2}{\xi^{2}} K_2(c|\xi|)\\
&=\dfrac{12}{\xi^4}-\dfrac{3c^2}{\xi^2}+6\left(\dfrac{1}{16}c^4\left (\dfrac{3}{2}-2\gamma\right)+\dfrac{1}{8}c^4\log(2)-\dfrac{1}{8}c^4\log(c)-\dfrac{1}{8}c^4\log(|\xi|)\right)+\cdots.
\end{align*}
Since we multiply by $(1-\cos(\xi))^2=\frac{\xi^4}{4}{\rm plus\ H.O.T.}$, $\xi\rightarrow 0$, we would have $\hat{\psi}(0)=\frac{12}{4}=3$. Therefore, in this case, we should divide by 3. 
After dividing by 3, it can be shown that 
$$\hat{\psi} (0)=1, \qquad  \dfrac{d\, \hat{\psi}}{d \xi}(0)=0, \qquad \dfrac{d^{k}\hat{\psi}}{d \xi^k}(2\pi j)=0, \: \forall j \in \mathbb{Z}\setminus \{0\}, \: k=0, 1.$$
As $\dfrac{d^2\, \hat{\psi}(\xi)}{d \xi^2}(0)\neq 0$ the Strang-Fix
condition for the second derivative is not satisfied at $\xi=0$. Therefore
we will have $\mathbb{P}_1$ reproduction at most.
\begin{theorem} Let $\phi$ be the generalised multiquadric radial basis
  function. With $n=1$ and the Fourier coefficients of  $\dfrac{1}{3}(2-2\cos
  x)^{2}$ being the $\lambda_j$, the quasi-interpolation 
  $$ \sum_{k\in\mathbb{Z}} p(k)\psi(\cdot-k)\equiv p$$
  is exact for linear polynomials $p$.
  \end{theorem}
The approximation error will have by  Theorem \ref{SF} with  $m=1,
\ell=3, n=1$ the asymptotic order $O(h^2)$.
\begin{remark}
As 
$$
(1-\cos\xi)^2=\frac{1}{4}e^{-i2\xi}-e^{-i\xi}+\frac{3}{2}-e^{i\xi}+\frac{1}{4} e^{i2\xi},\qquad -\pi\leq \xi\leq\pi,
$$
the kernel $\psi$ is given by the finite sum
$$
\psi(x)=\frac{1}{4}\varphi(|x+2|)-\varphi(|x+1|)+\frac{3}{2}\varphi(|x|)-\varphi(|x-1|)+\frac{1}{4} \varphi(|x-2|),\;x\in\RR.
$$
We can also use Mathematica  to check that this finite sum satisfies
$$|\psi (x)|=O((1+|x|)^{-5}),\quad |x|\to \infty.$$
\end{remark}
\subsection{Improving the reproduction}
Starting with the above explicit construction, we can improve the above scheme to $\mathbb{P}_3$-reproduction changing $(1-\cos (\xi))^2$ by a suitable trigonometric polynomial $P(\xi)$. We consider
\begin{itemize}
\item as the radial basis function we take the generalised multiquadric function defined by  $\varphi(r)=(r^2+ c^2)^{3/2}$,
\item as the finite sum of exponentials in order to resolve the
  generalised multiquadrics' singularity  $P(\xi)=\displaystyle\sum_{k=-N}^{N} \mu_k e^{ik\xi}$,
\item and therefore, finally, as Fourier transform $\hat{\psi}(\xi)=P(\xi) \hat{\varphi}(|\xi|)$.
\end{itemize}
We adjust the coefficients $\mu_k$ of $P(\xi)$ in such a way that the
expansion of $\hat{\psi}(\xi)$ around $\xi=0$ is
$1+O\left(\xi^4\log(|\xi|\right)$. To do that, we impose that, in that expansion,
the coefficient of $\xi^0$ is 1 and the coefficients of
$\xi^{-4},\xi^{-3},\xi^{-2},\xi^{-1},\xi,\xi^2,\xi^3$ vanish (even the
vanishing of the coefficient of $\xi^4$ may be added, but not the one
of $\xi^4 \log (|\xi|)$, which is not compatible with the previous
conditions due to the log-term). The system to be solved collecting all of these
requirements is
\begin{equation*}
\begin{split}
\sum_{k=-N}^N \mu_k=\sum_{k=-N}^N k\mu_k=\sum_{k=-N}^N k^2\mu_k= \sum_{k=-N}^N k^3\mu_k=0, \;  \sum_{k=-N}^N k^4\mu_k=2, \quad \sum_{k=-N}^N k^5\mu_k=0, \; \\ \sum_{k=-N}^N k^6\mu_k=-15c^2, \; \sum_{k=-N}^N k^7\mu_k=0, \; \sum_{k=-N}^N k^8\mu_k=\frac{105}{2}c^4(1+4\gamma-4\log 2+4\log c).
\end{split}
\end{equation*}
Its solution for $N=4$ is:
\begin{equation*}
\begin{split}
\mu_{-4}=\mu_4=&\frac{1}{11520}(28+60c^2+15c^4+60c^4\gamma-60c^4\log 2+60c^4\log c), \\ \mu_{-3}=\mu_3=&\frac{1}{480}(-16-30c^2-5c^4-20c^4\gamma+20c^4\log 2 -20c^4\log c),\\
\mu_{-2}=\mu_2=&\frac{1}{2880} (676+780c^2+105c^4+420c^4\gamma-420c^4\log 2+420c^4\log c), \\ 
\mu_{-1}=\mu_1=&\frac{1}{1440}(-976-870c^2-105c^4-420c^4\gamma+420c^4\log 2-420c^4\log c),\\
\mu_0=&\frac{1}{384} (364+300c^2+35c^4+140c^4\gamma-140c^4\log 2+140 c^4 \log c).
\end{split}
\end{equation*}
By inserting these coefficients into $P(\xi)$ we obtain that the expansion of $\hat{\psi}(\xi)$ around $\xi=0$ is
$$
1-\frac{c^4}{16}\xi^4\log |\xi|+\cdots.
$$
Therefore, by applying the (inverse) Fourier transform (see, e.g., \cite{jon}), we obtain a {\it decay of order -5 for $\psi$:
$$|\psi (x)|=O((1+|x|)^{-5}),\quad |x|\to \infty.$$
}
On the other hand, for this choice of $P(\xi)$, the (distributional) Fourier transform of $\psi$,
$$\hat{\psi}(\xi)=(2\pi)^{1/2} \frac{2^{5/2}}{\Gamma (-3/2)} \frac{c^2}{\xi^{2}} K_2(c|\xi|) P(\xi),\qquad \xi\in\RR,$$
satisfies the conditions
\begin{equation*}
\begin{split}
\hat{\psi}(0)=1,\quad  \frac{d^k \hat{\psi}}{d\xi^k}(2\pi j)=0, \, \forall j\in\mathbb{Z},\, k=1,2,3.
\end{split}
\end{equation*}
Therefore, according to the Strang-Fix conditions,  $\mathbb{P}_3$-reproduction is reached by the quasi-interpolant. As $\frac{d^4 \hat{\psi}}{d\xi^4}(0)=\infty$ it is not possible $\mathbb{P}_4$-reproduction in this case.
\begin{theorem} Let $\phi$ be the generalised multiquadric
  function. With $n=1$ and the explicit $\lambda_j$ as above for
  general parameters $c$, the quasi-interpolation 
  $$ \sum_{k\in\mathbb{Z}} p(k)\psi(\cdot-k)\equiv p$$
  is exact for cubic polynomials $p$.
  \end{theorem}
Putting all together, the {\it approximation order of the quasi-interpolant is $h^4 |\log h|$} (Theorem \ref{SF} with $m=3, \ell=1, n=1$).

\begin{remark}
The kernel $\psi$ is given by the finite sum
$$
\psi(x)=\sum_{k=-4}^4 \mu_k \varphi (|x-k|),\qquad x\in\RR,
$$
where $\mu_k$, $k=-4,\ldots,4$ are the ones as above. Using Mathematica we are also able to check that this finite sum satisfies
$$|\psi (x)|=O((1+|x|)^{-5}),\quad |x|\to \infty.$$
\end{remark}

\subsection{Scheme with the cubic B-spline}
The basic ideas of forming quasi-interpolants, quasi-Lagrange
functions and providing polynomial precision come of course from
B-spline quasi-interpolation. This one has the same behaviour than the function of the previous section at $r=0,$ i.e., we set
\begin{itemize}
\item as a radial basis function  $\varphi(r)=r^3$ and therefore we
  have the Fourier transform $ \hat{\varphi}(r)$ to be a constant multiple of $r^{-4}$,
\item the quasi-Lagrange functions $\psi$ are of course the compactly
  supported normalised
  B-splines that have well-known analytic Fourier transforms $\hat{N_k}(\xi) = \left( \hat{N_0}(\xi) \right)^{k+1}$ being $\hat{N_0}(\xi)=\dfrac{1-e^{-i \xi}}{i \xi}$,
\item  and our well-known coefficients $\lambda_j$  are the (finitely supported) Fourier coefficients of $(1-\cos (\xi))^{2}$. 
\item the decay of $\psi(x)$ is clear, we notice that it is faster
  than $O((1+ |x|)^{-k})$ for any $k\in \mathbb{N}$.
\end{itemize}
Studying this setup using Theorem \ref{SF} with  $\ell  > 1, n=1$ gives
\begin{itemize}
\item as far as polynomial precision is concerned, 
  $\mathbb{P}_1$-reproduction. The same case as in the first
  part of Subsection \ref{2s3}. By the way, we could improve the 
  polynomial exactness to $\mathbb{P}_3$ by adding a suitable polynomial $P(x)$.
\item Finally, the approximation error is 
$O(h^{m+1})$, the exponent being $m=1$ or $m=3$ depending on
polynomial reproduction order in our particular cases. 
\end{itemize}

\section{A further construction in the Fourier domain}\label{sec4}

Looking closely again at the Strang-Fix conditions we notice that the second and third set of conditions can be easily satisfied starting the construction in the Fourier domain. When we start in the Fourier domain the more important part is to ensure that the first condition is satisfied. In order to deduce the decay conditions of the function from the Fourier transform  we establish the following result. 

\begin{theorem}\label{th1}
Let $f$ be a symmetric real valued  function satisfying
\begin{enumerate}
\item $f\in C^{\underline n-1}(\mathbb{R})$ with all derivatives absolutely
  integrable and
\item $f^{(\underline n)}$ is a locally absolutely continuous function
and it has bounded variation on $\mathbb{R},$
\end{enumerate}
then $|\hat{f}(\omega)|= o(|\omega|^{-\underline n-1})$, $\omega \to \pm
\infty$.
\end{theorem}
\begin{proof}
  It follows immediately from the Riemann-Lebesgue Lemma and by integration by parts that
  $|\hat f(\omega)|=o(\omega^{-\underline n})$, $\omega\to\pm\infty$.

   Also, $f^{(\underline n)}$ has a
  well-defined Fourier transform with bounded $L^1$-norm (see  the
  first display  in \cite[Theorem 1]{lif}) and it even has an
  integrable derivative, so we have once more by the Riemann-Lebesgue Lemma
  that even $|\hat f(\omega)|=o(|\omega^{-\underline n-1}|)$. Notice
  for the argument of partial differentiation where we integrate the
  exponential and differentiate the other factor in the integrand,
  integrability is sufficient. This is what we wanted to show.
  \end{proof}

\subsection{A first example}
By Theorem \ref{th1}, it will be useful to have a
function $\psi$ with  Fourier transform of type
$\hat{\psi}(\xi)=\dfrac{Q(\xi)}{|\xi|^3},$ with $Q(\xi)$ such that
$\hat{\psi}(\xi)$ satisfies the following conditions:
\begin{enumerate}
\item It belongs to $C^2(\mathbb{R})$ and $\hat{\psi}'''(\xi)$ is a
  piecewise continuous function, and it has bounded variation in
  $\mathbb{R}$. Moreover $\hat{\psi}(\xi)$ and its derivatives up to  the third order have to be integrable.
\item The Strang-Fix conditions must be satisfied in order to get $\mathbb{P}_2-$reproduction which requires that $Q(\xi)$ must have a zero of order 3 at the origin.
\end{enumerate}
We think that the point is that we must avoid jumps in the Fourier
transform. More in detail,  if $\hat{f}^{(\underline n)}(\omega)$ has a jump then
it will have a behaviour like a multiple and shift of the Heaviside
function. So, we will obtain for $f(x)$ an asymptotic decay of  $O(|x|^{-\underline{n}-1})$ at maximum (although
we also need some other conditions as absolute integrability and
bounded variation of the function and its derivatives). The reason is
that the derivative picks up a Dirac delta, which has a Fourier
transform of constant modulus.


With the above conditions, we have found the example $$\hat{\psi}(\xi)= e^{-a \xi^4+ \frac{\xi^2}{2}}\left |\dfrac{\sin \xi}{\xi}\right|^3, \quad a>0.$$ 
The function in the Fourier domain is plotted in Figure 1. 
This function satisfies all our requirements. Moreover, as the third derivative of  $\hat{\psi}(\xi)$ has a jump at $\xi=\pi$, visualised in Figure 2, we would obtain an asymptotic  decay of $(1+|y|)^{-4}$.
In fact,  computations with Mathematica indicate that this is so. We have the following properties: 
\begin{enumerate}
\item  Decay of $\psi(x)$ turns out to be   $\psi(x)=o((1+|x|)^{-4})$: to find that, we apply Theorem \ref{th1} with $\underline n=3$.
\item Polynomial precision: $\mathbb{P}_2-$reproduction.
\item Approximation error: $O(h^3 |\log h|)$, (use Theorem \ref{SF} with  $m=2, \ell=1, n=1$).
\end{enumerate}
 \begin{figure}[h!]\label{invGau}
\includegraphics[width=0.4\textwidth]{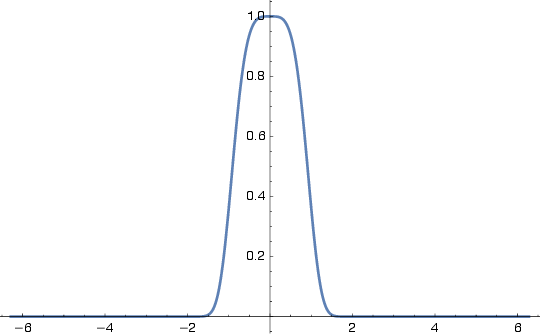}
\hskip1cm
\includegraphics[width=0.4\textwidth]{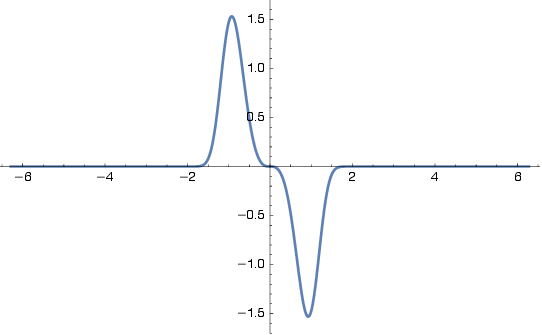}
\caption{Graphs of $\hat{\psi}(\xi)$ (left) and $\hat{\psi'}(\xi)$ (right)}
\end{figure}

\begin{figure}\label{Jump}
\includegraphics[width=0.3\textwidth]{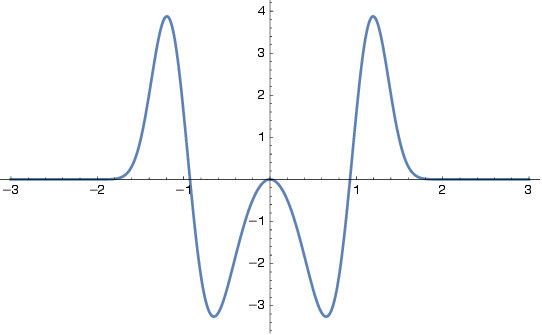}
\hskip0.1cm
\includegraphics[width=0.3\textwidth]{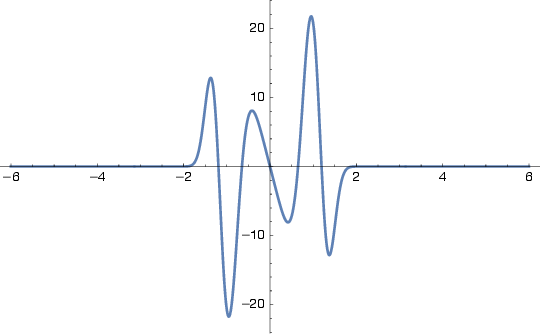}
\hskip0.1cm
\includegraphics[width=0.3\textwidth]{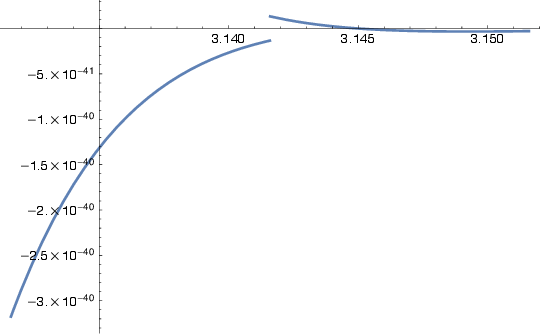}
\caption{Graphs of $\hat{\psi''}(x)$ (left) and $\hat{\psi'''}(x)$ (center) and  $\hat{\psi'''}(x)$ (right) around $x=\pi$.}
\end{figure}
 
 \subsection{Generalisation of the construction}
We now give a general description of a method to derive similar Quasi-Lagrange functions with higher order polynomial reproduction. 
We look for a $\psi$ with Fourier transform of type
\begin{equation}\label{Construct}
\hat{\psi}(\xi)=\dfrac{Q(\xi)\vert \sin (\xi)\vert ^{m+1}}{|\xi|^{m+1}},
\end{equation} with $Q(\xi)=e^{p(\xi)}$, for a polynomial $p$. Further
$\hat{\psi}(\xi)$ should satisfy the following conditions:
\begin{enumerate}
\item It belongs to $C^m(\mathbb{R})$. Moreover, $\hat{\psi}^{m+1}(\xi)$ is a
  piecewise continuous function and  has bounded variation on
  $\mathbb{R}$. Furthermore, $\hat{\psi}(\xi)$ and its derivatives up to order $m$ have to be integrable.
\item The Strang-Fix conditions cited in Theorem  \ref{SF} must be satisfied in order to get $\mathbb{P}_m-$reproduction.
\end{enumerate}
The properties can be translated into conditions on $Q$, from \eqref{Construct}.
Since for $m$ even
$$g(\xi)=\vert \sin (\xi)\vert^{m+1}=\sin(\xi)^{m+1}  \operatorname{sign}(\sin(\xi))$$
 is in $C^{\infty}(\RR \setminus \pi \ZZ)$, it is further in  $C^{m}(\RR)$  since $g^{(j)}(k\pi)=0$ for all $k\in \ZZ$ and $j\leq m$.

Combining this with the higher order product rules and continuity of $\vert \xi \vert^{-m-1}$ outside zero, $\hat{\psi}(\xi)$ satisfies the third condition of Theorem \ref{SF} if $Q(\xi)$ is $C^ {\infty}(\RR)$. 

Part 2 of the second condition  is satisfied if $Q(0)=1$, which is equivalent to assuming $p(0)=0$. For the first part of the second condition of Theorem \ref{SF}, we need to compute the derivative of $\hat{\psi}(\xi)$ near zero we find that for $\vert \xi \vert< \pi$:
\begin{equation*}
\begin{split}
\hat{\psi}'(\xi)=&Q'(\xi) \dfrac{\vert \sin(\xi)\vert^{m+1}}{|\xi|^{m+1}}
+Q(\xi) \dfrac{(m+1)\vert \sin(\xi)\vert^{m} \cos(\xi) \operatorname{sign}(\sin(\xi))}{|\xi|^{m+1}}\\ &
+Q(\xi) \dfrac{\vert \sin(\xi)\vert^{m+1} (-m-1)
  \operatorname{sign}(\xi)}{|\xi|^{m+2}}\\
=&Q(\xi) \dfrac{p'(\xi)  \sin(\xi)^{m+1}\xi
+ (m+1)\sin(\xi)^{m} \cos(\xi) \xi
+\sin(\xi)^{m+1} (-m-1)}{\xi^{m+2}},
\end{split}
\end{equation*}
where we used that for $\vert \xi \vert < \pi$, $\operatorname{sign}(\sin \xi)=\operatorname{sign}(\xi)$.
The condition  $\hat{\psi}'(0)=0$ is therefore satisfied if
$$p'(\xi) \sin(\xi)^{m+1}\xi
+ (m+1)\sin(\xi)^{m} \left( \cos(\xi) \xi
-\sin(\xi) \right)= {O}(\xi^{m+3}),$$
taking into account that $\lim \frac{\sin(x)}{x}=1$ for $x\rightarrow 0$. It remains to show that
$$p'(\xi) \sin(\xi)\xi
+  (m+1)  \left( \cos(\xi) \xi
- \sin(\xi)\right)= {O}(\xi^{3}).$$
Using the Taylor series of the trigonometric functions we can show
$$p'(\xi) \left( \xi^2-\frac{\xi^4}{6}\right) + (m+1)  \left( -\frac{1}{3}\xi^3+\frac{1}{30}\xi^5 \right) + H.O.T. ={O}(\xi^{3}),$$
which is true if 
$$p'(\xi) = {O}(\xi).$$
 This condition can be satisfied for $\ell=1$ if $p(\xi)=\sum_{k=2}^{j} a_k \xi^k$ for any choice of coefficients $a_k$ with $a_j<0$.
 For higher order derivatives one needs to ensure that the Taylor expansion near zero of
\begin{equation}\label{conditionQ} Q(\xi)\sin(\xi)^{m+1}=\xi^{m+1}+O(\xi^{2m+2}).
\end{equation}
\begin{lemma}
Let $\hat{\psi}$ be of the form $\eqref{Construct}$ and $Q(\xi)=e^{p(\xi)}$ satisfies \eqref{conditionQ} for an even $m$. Then  the quasi-interpolation 
  $$ \sum_{k\in\mathbb{Z}} p(k)\psi(\cdot-k)\equiv p$$
  is exact for polynomials $p$ of degree $m$.
\end{lemma}

\begin{remark}
\begin{enumerate}
\item For the case of $m=2$ this gives
 \begin{align*}
 Q(\xi)\sin(\xi)^{m+1}&=e^{p(\xi)}\left(\frac{\xi}{1}-\frac{\xi^3}{6}+...\right)^3\\
 &=\left( 1+ p'(0)\xi+\left(p''(0)+p'(0)^2\right)\frac{\xi^2}{2}+\cdots \right)\left(\frac{\xi}{1}-\frac{\xi^3}{6}+\cdots\right)^3\\
 &=\left( 1+ p''(0)\frac{\xi^2}{2}\right)\left(\xi^3-\frac{\xi^5}{2}+\cdots \right)\\
\end{align*}
which is equal to $\xi^{3}+O(\xi^{6})$ if we choose $p(\xi)=-\xi^4+\frac{\xi^2}{2}$.
\item In order to compute the function $\psi$ it is helpful to note that the inverse Fourier transform of $e^{-x^{4}}$ has been computed in \cite{Boyd2014} and a series representation of the Fourier transform of $e^{-\Vert x\Vert^{\beta}}$ is given in the Appendix A.2 to also compute higher order polynomial reproduction properties.
\end{enumerate}
\end{remark}

%


\section{Summary of the results}
In Section \ref{sect2} we compared two options of constructing quasi-interpolants using thin-plate spline and an infinite number of linear combinations to cardinal interpolation. The results are summarised in Table \ref{TPS}.

In Section \ref{sect3} we studied multiquadric and polyharmonic-spline
based quasi-interpolants. The result for the use of these basis
functions with and without forming adequate finite linear combinations
are displayed in Table \ref{MQ}. Two alternative constructions can be found in Section \ref{sec4} - where a new radial basis function characterised in the Fourier
domain which gives good  approximation properties is introduced - (see Table
\ref{New}). For the general construction described in this section the approximation order can even be increased to order $h^{m+1} |\log h|$.
\begin{table}[h]
\centering
\begin{tabular}{ |c|c|c|c|c| } 
 \hline
 Method & $\lambda$s & Reproduction & Decay & Approximation order \\ \hline\hline
 Cardinal  & Infinite & $\mathbb{P}_2$ & $O((1+|x|)^{-4})$ & $O(h^3 |\log h|)$\\ 
 Subsection \ref{subsect22}   & Infinite & $\mathbb{P}_1$ & $O((1+|x|)^{-4})$ & $O(h^2)$\\ 
 Subsection \ref{subsect23} & Infinite & $\mathbb{P}_2$ & $O((1+|x|)^{-2})$ & $O(h |\log h|)$\\ 
 \hline
\end{tabular}
\caption{$\varphi(r)= r^2\log r$}\label{TPS}
\end{table}

\begin{table}[h]
\centering
\begin{tabular}{ |c|c|c|c|c| } 
 \hline
 RBF $\phi(r)$& $\lambda$'s & Reproduction & Decay & Approximation order \\ \hline\hline
  $(r^2+ c^2)^{3/2}$ & Finite & $\mathbb{P}_1$ & $O((1+|x|)^{-5})$ & $O(h^2)$\\ 
$(r^2+ c^2)^{3/2}$ & Finite$^\ast$  & $\mathbb{P}_3$ &  $O((1+|x|)^{-5})$ & $O(h^4 |\log h|)$\\ 
Cubic B-spline & Finite & $\mathbb{P}_1$ & $O((1+|x|)^{-k}), \, \forall k\in\mathbb{N}$ & $O(h^2)$\\ 
Cubic B-spline & Finite$^\ast$ & $\mathbb{P}_3$ & $O((1+|x|)^{-k}), \, \forall k\in\mathbb{N}$ & $O(h^4)$\\ 
 \hline
\end{tabular}
\caption{The asterisk means that we included a suitable polynomial.}\label{MQ}
\end{table}

\begin{table}[h]
\centering
\begin{tabular}{ |c|c|c|c|c| } 
 \hline
 Function & $\lambda$'s & Reproduction & Decay & Approximation order \\ \hline\hline
$\hat{\psi}(\xi)$ &Infinite  & $\mathbb{P}_2$ & $O((1+|x|)^{-4})$ & $O(h^3 |\log h|)$\\ 
 \hline
\end{tabular}
\caption{The  function $\hat{\psi}(\xi)= e^{-a x^4+ \frac{x^2}{2}}\left |\dfrac{\sin \xi}{\xi}\right|^3, \quad a>0$.}\label{New}
\end{table}

\clearpage
\section*{Appendix}

\subsection*{A.1 Computation of the inverse Fourier transform of $\vert \sin(\xi)\vert$}

\begin{lemma}\label{sin}
The inverse Fourier transform of $|sin(x)|$ is
  $${\cal F}^{-1}|\sin|(x)=\frac1{\pi}\times\frac{1+\exp(ix\pi)}{1-x^2}\times{\cal D}_2(x).$$
  Here ${\cal D}_2$ is the Dirac comb
  $${\cal D}_2=\sum_{k=-\infty}^\infty \delta(\cdot-2k).$$
\end{lemma} 
\begin{proof}
In Section \ref{sect2} we used that if
 $g(\xi)=\frac12\sin|\xi|+\frac12\sin|\xi-\pi|$, then
  $$|\sin \xi| = \sum_{k=-\infty}^\infty g(\xi-k\pi),\qquad
  \xi\in\mathbb{R}.$$
  In order to give the inverse Fourier transform of $|\sin \xi|$ we note first that the generalised inverse Fourier transform ${\cal F}^{-1}$ of $\sin|\xi|$ is
   \begin{align*}
   {\cal F}^{-1}\sin|\cdot|(x) &=\frac{1}{2\pi}\int_{-\infty}^{\infty} \sin|\xi| \exp(i\xi x) \,d\xi 
   =\frac{1}{\pi}\int_{0}^{\infty}\sin \xi \cos(\xi)\,d\xi \\
&= \frac{1}{\pi}\lim_{\epsilon\to0_+}\int_0^\infty\exp(-\epsilon \xi)\sin \xi\cos \xi x\,d\xi \\
&= \frac{1}{\pi}\sqrt{\frac{\pi}{2}}\lim_{\epsilon\to0_+}\sqrt{\frac{2}{\pi}}\int_0^\infty\exp(-\epsilon \xi)\sin \xi\cos \xi x\,d\xi \\
&= \frac{1}{\pi}\sqrt{\frac{\pi}{2}}\lim_{\epsilon\to0_+}{\cal F}_{cos}\left(\exp(-\epsilon \xi)\sin \xi\right) (x)\\
&=\frac{1}{\sqrt{2\pi}}\lim_{\epsilon\to0_+}\frac{1}{\sqrt{2\pi}} \left( \frac{1+x}{\epsilon^2+(1+x)^{2}} +\frac{1-x}{\epsilon^2+(1-x)^2}\right)\\
  &=\frac{1}{2\pi} \frac{2}{1-x^2}=\frac{1}{\pi}\frac1{1-x^2}.
  \end{align*}
  where we used the cosine transform, ${\cal F}_{cos}$, given in \cite{GR},  ($17.34.22^{7}$). This gives
  $${\cal F}^{-1}\sin|\cdot|(x)=\frac1\pi\frac1{1-x^2}.$$
  Therefore, the generalised inverse Fourier transform of $g $ is
  $${\cal F}^{-1}g(x)=\frac1{2\pi}\biggl(\frac1{1-x^2}+\frac{\exp(ix\pi)}{1-x^2}\biggr).$$
  This gives according to \cite{SW}
  $${\cal F}^{-1}|\sin|(x)=\frac1{\pi}\times\frac{1+\exp(ix\pi)}{1-x^2}\times{\cal D}_2(x).$$
      First, we have
      $$ |\sin (x)|=\sum_{k=-\infty}^\infty g(x-k\pi)=\left(\sum_{k=-\infty}^{\infty} \delta(\cdot- k\pi)\ast g\right)(x)$$
Now, from $\mathcal{F}^{-1}f(x) = \frac{1}{(2\pi)^n} (\mathcal{F}f)(-x)$ and 
$\mathcal{F}^{-1}(f*g)(x)= (2\pi)^n \, \mathcal{F}^{-1}f(x) \cdot  \mathcal{F}^{-1}g(x)$ 
it follows:
  \begin{align*}
 {\cal F}^{-1}  |\sin (\cdot)|(x)=&2\pi{\cal F}^{-1}  \left(\sum_{k=-\infty}^{\infty} \delta(\cdot- k\pi)\right)(x) \times {\cal F}^{-1} g(x)
 \\=& 2\pi \biggl(\frac{1}{\pi} \left(\sum_{k=-\infty}^{\infty} \delta(\cdot- 2k)\right)(x)  \biggr) \times \biggl(  \frac1{2\pi}\biggl(\frac1{1-x^2}+\frac{\exp(ix\pi)}{1-x^2}\biggr) \biggr) 
 \\=& {\cal D}_2(x) \times \frac1{\pi}\biggl(\frac{1+\exp(ix\pi) }{1-x^2}\biggr).
  \end{align*}
  \end{proof}
  
  \subsection*{A.2 Inverse Fourier transform of $\exp(-\Vert x\Vert^{\beta})$}
  We now want to investigate the class of inverse $n$-dimensional Fourier transforms of the functions 
\begin{equation}Q(\xi)=e^{-\Vert \xi\Vert ^{\beta}},\end{equation}
which are integrable for $\beta>0$.
The presented results are based on the results given in the thesis of one of the authors 
\cite[Chapter 3.3]{Jaeger2018}.  We start by gathering informations about the special choices of $\beta$ which  have already been considered.
  
\begin{itemize}
\item $\beta=1$: In this case the function is  $$Q( \xi )=e^{-\Vert \xi \Vert},$$ which is the Poisson kernel. Its Fourier transform is  $$\mathcal{F}^{-1}{Q}(\xi)=\frac{1}{2\pi} \Gamma\left(\frac{n}{2}+\frac{1}{2}\right)\frac{1}{(1+\Vert \xi \Vert^2)^{\frac{n}{2}+\frac{1}{2}}},$$ which is a special case of the generalised inverse multiquadric,
 \mbox{$\phi(r)=(1+r^2)^{\alpha/2}$}, with $\alpha=-n-1$,
\item $\beta=2$: The function is the Gaussian basis function $Q(\xi)=e^{-\Vert \xi \Vert ^2}$, which has the inverse Fourier transform $\mathcal{F}^{-1}{Q}(\xi)=(1/4\pi)^{n/2}e^{-\Vert \xi \Vert^2/4}$ which is also a Gaussian basis function,
\item $\beta=2N$:  The function is ${Q}(\xi)=e^{-\Vert \xi \Vert^{2N}}$; its Fourier transform was considered, for the case $n=1$, in \cite{Boyd2014}.  The Fourier transforms of $Q(\xi)=e^{-A\vert \xi \vert ^{2N}}$ have therein been approximated  without giving a representation different from the obvious integral description.  For the special case $\beta=4$ the resulting radial basis function is called the inverse quartic Gaussian ($\beta=4$). A series representation has been computed using Matlab by Boyd in \cite{Boyd2013} and takes the form
\begin{multline}\label{eq:Boyding} \mathcal{F}^{-1}{Q}(\xi)=\frac{1}{2^{3/2}}\sum_{k=0}^{\infty} \frac{\Gamma(1/2)}{\Gamma(1/2+N)\Gamma(3/4+k)}\frac{\left(\frac{\vert \xi\vert}{4}\right)^{4k}}{k!}\\
-\frac{1}{8\pi}\Gamma(3/4)\vert \xi\vert^2 \sum_{k=0}^{\infty} \frac{\Gamma(5/4)\Gamma(3/2)}{\Gamma(3/2+k)\Gamma(5/4+k)}\frac{\left(\frac{\vert \xi\vert}{4}\right)^{4k}}{k!}.
\end{multline}
\end{itemize}

  We now give a representation of the inverse Fourier transform of $Q(\xi)=e^{-\Vert \xi\Vert^{\beta}}$.
   
We focus on the case $\beta > 1$ using the series representation of the Bessel function. 
However, to be able to compute the Fourier transform we need to prove this additional lemma first.
\begin{lemma}\label{leInvGaus2}
The series 
\begin{equation}\sum_{k\geq 0}(-1)^ka^{2k}\frac{\Gamma\left( \frac{n+2k}{\beta}\right)}{\Gamma(k+1)\Gamma(k+ \frac{n}{2})}, \quad a\in \mathbb{R},
\end{equation}
is absolutely convergent for every $\beta >1.$
\end{lemma}
\begin{proof}
We  are going to prove that  by applying the root test to the series the resulting limit is 0. First, we have
\begin{equation}\label{lim}
0 \leq \lim_{k\to \infty} \left | (-1)^ka^{2k}\frac{\Gamma\left( \frac{n+2k}{\beta}\right)}{\Gamma(k+1)\Gamma(k+ \frac{n}{2})}\right |^{\frac1k}
=\lim_{k\to \infty} \left (a^{2k}\frac{\Gamma\left( \frac{n+2k}{\beta}\right)}{\Gamma(k+1)\Gamma(k+ \frac{n}{2})}\right )^{\frac1k}.
\end{equation}
Applying formula 8.327, $1^{\ast}$  of \cite{GR} we have that for large $k$
\begin{align*}
 \left( \frac{ a^{2k}\Gamma\left( \frac{n+2k}{\beta}\right)}{\Gamma(k+1)\Gamma(k+  \frac{n}{2} )}\right)^{\frac1k}
 &\leq  \left( \frac{a^{2k}\sqrt{2\pi} \sqrt{\frac{n+2k}{\beta}-1} \left(\frac{n+2k}{\beta}-1\right)^{\frac{n+2k}{\beta}-1}  e^k e^{k+ \frac{n}{2}  -1}c}{ \sqrt{2\pi} \sqrt{k} k^{k}
 \sqrt{2\pi} \sqrt{k+  \frac{n}{2}-1} \left(k+ \frac{n}{2}-1\right)^{k+  \frac{n}{2} -1} 
 e^{\frac{n+2k}{\beta}-1}}\right)^{\frac1k} \\
 &= \left (\frac{c e^{n(\frac12- \frac{1}{\beta})} \left(\frac{n+2k}{\beta}-1\right)^{\frac{n}{\beta}-\frac12}
 }{ \sqrt{k}
 \sqrt{2\pi} \left(k+  \frac{n}{2}-1\right)^{ \frac{n}{2}-\frac12}}\times
\dfrac{ a^{2k}e^{2k(1- \frac{1}{\beta})} 
 \left(\frac{n+2k}{\beta}-1\right)^{\frac{2k}{\beta}}
 }{k^k (k+  \frac{n}{2}-1)^{k}}\right)^{\frac1k} \\
 &=\left (\frac{c e^{n(\frac12- \frac{1}{\beta})} \left(\frac{n+2k}{\beta}-1\right)^{\frac{n}{\beta}-\frac12}
 }{ \sqrt{k}
 \sqrt{2\pi} \left(k+  \frac{n}{2}-1\right)^{ \frac{n}{2}-\frac12}}\right)^{\frac1k} \times\left (
\dfrac{ a^{2}e^{2(1- \frac{1}{\beta})} 
 \left(\frac{n+2k}{\beta}-1\right)^{\frac{2}{\beta}}
 }{k (k+  \frac{n}{2}-1)}\right)
 \end{align*}
 where $c$  is a constant that comes from the remainder of the series in  the numerator. The last expression, for large $k$ involves a  (first) fraction that tends to one and a second fraction that is of order  $O(k^{\frac{2}{\beta}-2}).$  As $\beta>1,$ the absolute convergence of the series for any $a\in \mathbb{R}$ follows.
\end{proof} 

\begin{lemma}\label{thm:invGau2}
The inverse Fourier transform of $Q(x )=e^{-\Vert x\Vert^{\beta}}$, $x\in \R^n$, $ \beta > 1$  is
\begin{equation}\mathcal{F}^{-1}Q(\xi)=\frac{2^{1-n}}{\pi^{n/2}}  \frac{1}{\beta}\sum_{k=0}^{\infty}\frac{(-1)^k\left(\frac{\Vert  \xi \Vert}{2}\right)^{2k}}{k!\Gamma(k+\frac{n}{2})}\Gamma\left(\frac{n+2k}{\beta}\right).\end{equation}
\end{lemma}
\begin{proof}
We use the formula for Fourier transforms of radial functions to compute the inverse Fourier transform; this is applicable because $Q \in L^1(\R^n)$, for all $\beta>1$, and $n\in \N$. We then use the series representation of the Bessel function (\cite{abr} (9.1.10))
\begin{align*}
\mathcal{F}^{-1}Q(\xi)&=\frac{1}{\sqrt{2\pi}^{n}}\Vert  \xi \Vert^{-(\frac{n-2}{2})}\int_0^{\infty}e^{-t^{\beta}}t^{n/2} J_{\frac{n-2}{2}}(\Vert  \xi \Vert t)\,dt\\
&=\frac{1}{\sqrt{2\pi}^{n}}\Vert  \xi \Vert^{-(\frac{n-2}{2})}\int_0^{\infty}e^{-t^{\beta}}t^{n/2}\sum_{k=0}^{\infty} \frac{(-1)^k(\Vert  \xi \Vert t/2)^{2k+\frac{n}{2}-1}}{k!\Gamma(k+\frac{n}{2})}\,dt\\
&=\frac{1}{\sqrt{2\pi}^{n}}2^{-\frac{n}{2}+1}\int_0^{\infty}e^{-t^{\beta}}t^{n-1}\sum_{k=0}^{\infty} \frac{(-1)^k\left(\frac{\Vert  \xi \Vert t}{2}\right)^{2k}}{k!\Gamma(k+\frac{n}{2})}\,dt\\
&=\frac{1}{\sqrt{2\pi}^{n}}2^{-\frac{n}{2}+1}\underset{u\rightarrow \infty}{\lim}\int_0^{u}e^{-t^{\beta}}t^{n-1}\underset{n\rightarrow \infty}{\lim}\sum_{k=0}^{n} \frac{(-1)^k\left(\frac{\Vert  \xi \Vert t}{2}\right)^{2k}}{k!\Gamma(k+\frac{n}{2})}\,dt.\\
\end{align*} 
The sum inside the integrand can be bounded as follows:
$$\left \vert \sum_{k=0}^{n} \frac{\left(-1\right)^k\left(\frac{\Vert  \xi \Vert t}{2}\right)^{2k}}{k!\Gamma(k+\frac{n}{2})}\right \vert \leq\sum_{k=0}^{\infty}\left\vert \frac{\left(\frac{\Vert  \xi \Vert t}{2}\right)^{2k}}{k!\Gamma(k+\frac{n}{2})}\right \vert \leq e^{\frac14 (\Vert  \xi \Vert t)^{2}}+c,\quad c>0,$$ 
with $\frac{1}{\Gamma(k+n/2)}<1$, for $(k+n/2)>2$, which gives an integrable majorant on $[0,u].$ Thereby we get
\begin{align*}
\mathcal{F}^{-1}Q(\xi)&=\frac{1}{\sqrt{2\pi}^{n}}2^{-\frac{n}{2}+1}\underset{u\rightarrow \infty}{\lim}\sum_{k=0}^{\infty}\frac{\left(-1\right)^k\left(\frac{ \Vert  \xi \Vert}{2}\right)^{2k}}{k!\Gamma\left(k+\frac{n}{2}\right)}\int_0^{u}e^{-t^{\beta}}t^{n-1+2k}\,dt\\
&=\frac{1}{\sqrt{2\pi}^{n}}2^{-\frac{n}{2}+1}\underset{u\rightarrow \infty}{\lim}\sum_{k=0}^{\infty}\frac{\left(-1\right)^k\left(\frac{ \Vert  \xi \Vert}{2}\right)^{2k}}{k!\Gamma(k+\frac{n}{2})}\int_0^{u^{\beta}}e^{-z}z^{\frac{n+2k}{\beta}-1}\frac{1}{\beta}\, dz\\
&=\frac{1}{\sqrt{2\pi}^{n}}2^{-\frac{n}{2}+1}\frac{1}{\beta}\underset{u \rightarrow \infty}{\lim}\sum_{k=0}^{\infty}\frac{\left(-1\right)^k\left(\frac{ \Vert  \xi \Vert}{2}\right)^{2k}}{k!\Gamma\left(k+\frac{n}{2}\right)}\gamma\left(\frac{n+2k}{\beta},u^{\beta}\right),
\end{align*}
where we have used the expression 8.350.1 of \cite{GR} in the last equality.
Here $\gamma(\cdot,\cdot)$ is the incomplete $\Gamma$-function. We know that $\gamma\left(\frac{n+2k}{\beta},u^{\beta}\right)\leq \Gamma\left(\frac{n+2k}{\beta}\right)$ for all $\beta>1$ and applying Lemma \ref{leInvGaus2} we get a convergent majorant. So, we have 
\begin{align*}
\mathcal{F}^{-1}Q(\xi)&=\frac{2^{1-n}}{\pi^{n/2}} \frac{1}{\beta}\sum_{k=0}^{\infty}\frac{(-1)^k\left(\frac{ \Vert  \xi \Vert}{2}\right)^{2k}}{k!\Gamma\left(k+\frac{n}{2}\right)}\Gamma\left(\frac{n+2k}{\beta}\right).
\end{align*}
\end{proof}
The last series is absolutely convergent for $\beta>1$ and can be further simplified for many values of $\beta$ by applying the doubling or tripling formulas for the Gamma function.   For the application in Section \ref{sec4} we are specifically interested in the case $n=1$
and $Q(x)=\exp(-|x|^{\beta})$. In this case, and if $\beta=2N, N\in \mathbb{N}$  our formula simplifies to:
\begin{equation}\mathcal{F}^{-1}Q(\xi)=\frac{1}{2N\sqrt{\pi}}\sum_{k=0}^{\infty}\frac{(-1)^k\left(\frac{| \xi |}{2}\right)^{2k}}{k!\Gamma(k+\frac{1}{2})}\Gamma\left(\frac{1+2k}{2N}\right).
\end{equation}

\bibliographystyle{plain} 
\bibliography{bib}
\end{document}